\documentclass[11pt,letterpaper]{article}
\usepackage[T1]{fontenc}
\usepackage{lmodern}
\usepackage{amsmath,amssymb,amsthm}
\usepackage[textwidth=5.7in,top=1.06in,bottom=1.06in]{geometry}
\usepackage[expansion=false]{microtype}
\usepackage{etoolbox}
\usepackage[hidelinks]{hyperref}

\patchcmd{\thebibliography}{\section*{\refname}}{\section{\refname}}{}{}

\newtheorem{theorem}{Theorem}[section]
\newtheorem{lemma}[theorem]{Lemma}
\theoremstyle{definition}
\newtheorem{example}[theorem]{Example}
 \newtheorem{proposition}{Proposition}
  \newtheorem{conjecture}{Conjecture}

\newcommand{\calF}{\mathcal F}
\newcommand{\calH}{\mathcal H}
\newcommand{\calS}{\mathcal S}

\newcommand{\subsetsatmost}[2]{\binom{[#1]}{\le #2}}
\newcommand{\doi}[1]{\mbox{\href{https://doi.org/#1}{\nolinkurl{#1}}}}
\numberwithin{equation}{section}

\title{A complete intersection theorem for sets of bounded size}
\author{Jing Gu\thanks{Department of Mathematics, Changzhou University,
Changzhou 213164, China.}
\qquad
Biao Wu\thanks{School of Mathematics and Statistics, Hunan Normal University,
Changsha 410081, China.\newline
E-mail addresses: \href{mailto:wu@hunnu.edu.cn}{\texttt{wu@hunnu.edu.cn}} (Biao Wu), \href{mailto:gujing23@cczu.edu.cn}{\texttt{gujing23@cczu.edu.cn}} (Jing Gu).}}
\date{}
\hypersetup{
  pdftitle={A complete intersection theorem for sets of bounded size},
  pdfauthor={Jing Gu and Biao Wu}
}

\begin{document}
\maketitle
\begin{abstract}
Let \(1\le t\le k\le n\) and \(2k<n+t\). We prove that a largest
\(t\)-intersecting family of subsets of \([n]\) of size at most \(k\)
is obtained by taking all such subsets that contain at least \(t+r\)
points of a fixed \((t+2r)\)-element set, for some
\(0\le r\le\left\lfloor\frac{n-t}{2}\right\rfloor\).
This settles a conjecture of Ahlswede, Bey, Engel and Khachatrian from 2002.
\end{abstract}

\noindent\textbf{Keywords.} Extremal combinatorics, Intersecting families, Boolean lattice, Shift.

\smallskip
\noindent\textbf{2020 Mathematics Subject Classification.} 05D05.

\section{Introduction}

A family \(\calF\) of subsets of \([n]=\{1,\ldots,n\}\) is
\(t\)-intersecting if \(|A\cap B|\ge t\) for all \(A,B\in\calF\),
including \(A=B\). A \(1\)-intersecting family is simply called
\emph{intersecting}. Write \(\subsetsatmost nk\) and \(\binom{[n]}{k}\)
for the families of subsets of \([n]\) of size at most \(k\)
and of size \(k\), respectively.
For \(0\le r\le \left\lfloor\frac{n-t}{2}\right\rfloor\), let
\[
 \calS_r=\{A\subseteq[n]:|A\cap[t+2r]|\ge t+r\},
 \qquad
 \calS_r^{\le k}=\calS_r\cap\subsetsatmost nk.
\]
Any two members of \(\calS_r\) share at least \(t\) points in
\([t+2r]\). The family \(\calS_0\) consists of all sets containing
\([t]\) and is called a \(t\)-star.

The Erd\H{o}s-Ko-Rado Theorem stands as one of the cornerstone results in extremal combinatorics.

\begin{theorem}[Erd\H{o}s-Ko-Rado,~\cite{EKR61}]\label{EKR61}
Let $t,k,n$ be positive integers such that $t \leq k\leq n$. If $ \mathcal F \subseteq \binom{[n]}{k} $ is $t$-intersecting,  then for $n\geq n_0(k,t)$, \[|\mathcal{F}|\leq \binom{n-t}{k-t}. \]
\end{theorem}

Let \(N_0(k,t)\) be the smallest integer \(n_0\ge k\) such that
the bound in Theorem~\ref{EKR61} holds whenever \(n\ge n_0\).
For \(1\le t<k\), the sharp threshold is
\(N_0(k,t)=(t+1)(k-t+1)\). Erd\H{o}s, Ko and Rado~\cite{EKR61}
proved this for \(t=1\), Frankl~\cite{Frankl78} for \(t\ge15\),
and Wilson~\cite{Wilson84} for all \(t\).

Ahlswede and Khachatrian~\cite{AK97} proved the complete intersection theorem.

\begin{theorem}[Ahlswede-Khachatrian,~\cite{AK97}]\label{AK97}
Let $t,k,n$ be positive integers such that \(1\le t\le k\le n\). If $ \mathcal F \subseteq \binom{[n]}{k} $ is $t$-intersecting,  then
\[
 |\calF|\le
 \max_{0\le r\le \left\lfloor\frac{n-t}{2}\right\rfloor}\left|\calS_r\cap\binom{[n]}{k}\right|,
\]
with equality for some \(\calS_r\cap\binom{[n]}{k}\).
\end{theorem}

Katona~\cite{Katona64} proved the full Boolean lattice theorem.
\begin{theorem}[Katona,~\cite{Katona64}]
Let \(1\le t\le n\). If \(\calF\subseteq2^{[n]}\) is \(t\)-intersecting, then
\[
 |\calF|\le|\calS_{\left\lfloor\frac{n-t}{2}\right\rfloor}|=
 \begin{cases}
  \displaystyle\sum_{i=(n+t)/2}^{n}\binom ni,
      & n-t\text{ even},\\[6pt]
  \displaystyle 2\sum_{i=(n+t-1)/2}^{n-1}\binom{n-1}{i},
      & n-t\text{ odd}.
 \end{cases}
\]
When \(n-t\) is even, \(\calS_{\left\lfloor\frac{n-t}{2}\right\rfloor}\) consists of all sets of size
at least \((n+t)/2\). When \(n-t\) is odd, it consists of all
sets containing at least \((n+t-1)/2\) points of \([n-1]\).
\end{theorem}

For sets of size at most \(k\), Ahlswede, Bey, Engel and
Khachatrian proposed the following analogue of the uniform
complete intersection theorem.

\medskip

\begin{conjecture}[{Ahlswede-Bey-Engel-Khachatrian,~\cite{ABEK02}}]\label{conj-ABEK}
Let \(1\le t\le k\le n\) and \(2k<n+t\).
If $\calF\subseteq\subsetsatmost nk$ is $t$-intersecting, then
\[
|\calF| \le \max_{0\le r\le \left\lfloor\frac{n-t}{2}\right\rfloor}|\calS_r^{\le k}|.
\]
\end{conjecture}

For \(t=1\), the following bound holds for every \(k\).

\begin{proposition}[{Ahlswede-Bey-Engel-Khachatrian,~\cite{ABEK02}}]\label{prop:truncated-intersecting}
Let $1\le k\le n$. If
$
\mathcal F\subseteq \binom{[n]}{\le k}
$
is an intersecting family, then
\[
|\mathcal F|
\le |\calS_0^{\le k}|=
\sum_{j=1}^{k}\binom{n-1}{j-1}.
\]
\end{proposition}

For completeness, we give a short proof illustrating how
uniform bounds can be combined with complementary pairs.
By the Erd\H{o}s--Ko--Rado theorem,
\begin{equation}\label{eq:layer-ekr}
\left|\calF\cap\binom{[n]}j\right|
\le\binom{n-1}{j-1}
\qquad(1\le j\le k,\ 2j\le n).
\end{equation}
When \(2k\le n\), summing these bounds gives the result.

When \(2k>n\), put \(m=n-k\) and consider
\[
\mathcal B=\{A\subseteq[n]:m\le|A|\le k\}.
\]
Its members form complementary pairs
\(\{A,[n]\setminus A\}\).
Each pair contains at most one member of \(\calF\),
and exactly one member containing \(1\). Thus
\[
|\calF\cap\mathcal B|
\le\frac{|\mathcal B|}{2}
=\sum_{j=\max\{1,m\}}^{k}\binom{n-1}{j-1}.
\]
The remaining members of \(\calF\) have size below
\(m<n/2\), so \eqref{eq:layer-ekr} gives
\[
|\calF|
\le
\sum_{j=1}^{m-1}\binom{n-1}{j-1}
+\sum_{j=\max\{1,m\}}^{k}\binom{n-1}{j-1}
=
\sum_{j=1}^{k}\binom{n-1}{j-1},
\]
where an empty sum is zero.

Furthermore, Ahlswede, Bey, Engel and Khachatrian~\cite{ABEK02} proved the conjecture for fixed \(t\) and fixed
\(0<\varepsilon<1/2\) when \(k\le(1/2-\varepsilon)n\) and \(n\)
is sufficiently large.
They also proved that the maximum size of a \(t\)-intersecting
subfamily of \(\subsetsatmost nk\) is asymptotic to
\(\max_{0\le r\le \left\lfloor\frac{n-t}{2}\right\rfloor}|\calS_r^{\le k}|\) when
\(t=\tau n+o(n)\) and \(k=\kappa n+o(n)\), with
\(0<\tau<\kappa<(1+\tau)/2\).

Our main theorem establishes the conjecture throughout the range
\(2k<n+t\).

\begin{theorem}\label{thm:main}
Let \(1\le t\le k\le n\) and \(2k<n+t\).
If \(\calF\subseteq\subsetsatmost nk\) is \(t\)-intersecting, then
\[
 |\calF|\le\max_{0\le r\le \left\lfloor\frac{n-t}{2}\right\rfloor}|\calS_r^{\le k}|.
\]
The bound is attained by \(\calS_r^{\le k}\) for some $0\le r\le \left\lfloor\frac{n-t}{2}\right\rfloor$.
\end{theorem}

Summing the uniform bounds need not give the exact maximum,
as the following example shows.

\begin{example}\label{ex:levels}
Let \(n=7\), \(k=4\), and \(t=2\). Then
\[
|\calS_0^{\le4}|=16,\qquad
|\calS_1^{\le4}|=17,\qquad
|\calS_2^{\le4}|=15.
\]
Theorem~\ref{thm:main} therefore gives a maximum of \(17\).
The separate uniform maxima for sets of sizes \(2,3,4\)
are \(1,5,15\), respectively, whose sum is \(21\).
\end{example}

\section{Preliminaries}\label{sec:tools}

\subsection{Shifts}

The shifting method goes back to Erd\H{o}s, Ko and
Rado~\cite{EKR61}. Frankl~\cite{Frankl87} developed
the technique further and gave a systematic account of it.
See Wang~\cite{Wang26} for more recent developments.

Let \(1\le m\le n\). For \(A\subseteq[m]\) and \(1\le i<j\le m\),
let
\[
 s_{ij}(A)=
 \begin{cases}
  (A\setminus\{j\})\cup\{i\}, & j\in A,\ i\notin A,\\
  A, & \text{otherwise}.
 \end{cases}
\]
The \((i,j)\)-shift of \(\calH\subseteq 2^{[m]}\) replaces each
\(A\in\calH\) by \(s_{ij}(A)\) if \(s_{ij}(A)\notin\calH\).
The family is \emph{shifted} if \(s_{ij}(A)\in\calH\) for every
\(A\in\calH\) and every \(i<j\). Shifts preserve \(t\)-intersection
and the number of members of each size. Repeated shifts produce
a shifted family~\cite{Frankl87}.

For integers \(a\ge0\) and \(d\), put
\[
 B(a,d)=\sum_{0\le j\le d}\binom aj,
\]
where \(\binom aj=0\) for \(j<0\) or \(j>a\), and an empty sum
is zero. An extension of \(T\subseteq[m]\) is a set
\(A\subseteq[n]\) with \(|A|\le k\) and \(A\cap[m]=T\).
For \(|T|=s\), the number of extensions is \(B(n-m,k-s)\).

Assign each \(s\)-element subset of \([m]\) the weight
\begin{equation}\label{eq:trace-rank-weight}
 w_s=(2^m+1)B(n-m,k-s)+1,\qquad 0\le s\le m.
\end{equation}
If \(s>k\), the number of extensions is zero and \(w_s=1\).

 For \(\calH\subseteq 2^{[m]}\), define
\[
 w(\calH)=\sum_{A\in\calH}w_{|A|}.
\]
Shifts preserve this total weight, so a maximum-weight
\(t\)-intersecting family can be chosen shifted.
Since every \(w_s>0\), such a family is inclusion-maximal among
the \(t\)-intersecting families. It is also upward closed, since
any superset of a member meets every member in at least \(t\) points.

\subsection{Generating sets}

Ahlswede and Khachatrian introduced the generating set method in
their proof of the complete intersection theorem~\cite{AK97}
and adapted it for the complete nontrivial-intersection
theorem~\cite{AK96}. The method has been used to bound
products of cross-\(t\)-intersecting families~\cite{ZW25,HLWZ26,BJ26}
and sums of their sizes~\cite{LZ25}. Huang, Peng and Wang~\cite{HPW26}
used standard generating sets for lexicographic initial segments
in the pairwise cross-intersecting sum problem.

Let \(X\) be a finite set. The family generated by
\(\mathcal G\subseteq 2^X\) is
\[
 \langle\mathcal G\rangle_X
 =\{A\subseteq X:G\subseteq A\text{ for some }G\in\mathcal G\}.
\]
Every upward-closed family \(\calF\subseteq 2^X\) is generated
by its inclusion-minimal members. It is \(t\)-intersecting if
and only if these minimal members are \(t\)-intersecting.

The \emph{trace} of \(A\subseteq X\) on \(I\subseteq X\) is
\(A\cap I\). If \(\calF\) is upward closed and all its minimal
members are contained in \(I\), then
\[
 A\in\calF \quad\Longleftrightarrow\quad A\cap I\in\calF.
\]
Conversely, if this equivalence holds for every \(A\subseteq X\),
then each minimal member \(A\) equals \(A\cap I\) and hence lies
in \(I\). For example, \(\calS_r\) is generated by the
\((t+r)\)-subsets of \([t+2r]\).

For \(1\le m\le n\) and \(\calH\subseteq2^{[m]}\), the family
of all extensions of its members is
\[
\widehat{\calH}
=\{A\in\subsetsatmost nk:A\cap[m]\in\calH\}.
\]
If \(\calH\) is \(t\)-intersecting, then so is \(\widehat{\calH}\).
Each \(T\in\calH\) has \(B(n-m,k-|T|)\) extensions, and distinct
traces give disjoint families of extensions. Hence
\begin{equation}\label{eq:traceweight}
 |\widehat{\calH}|=\sum_{T\in\calH}B(n-m,k-|T|).
\end{equation}
With the weights in \eqref{eq:trace-rank-weight}, we obtain
\[
 w(\calH)=(2^m+1)|\widehat{\calH}|+|\calH|.
\]
Since \(|\calH|\le2^m\), maximizing \(w(\calH)\) first maximizes
\(|\widehat{\calH}|\) and, in case of a tie, \(|\calH|\).

Put \(f(q)=|\calS_q^{\le k}|\) for \(0\le q\le \left\lfloor\frac{n-t}{2}\right\rfloor\).

\begin{lemma}[{\cite{ABEK02}}]\label{lem:unimodal}
Let \(2\le t\le k\le n\) and \(2k<n+t\).
Then the sequence of \(f(q)\) for
\(0\le q\le \left\lfloor\frac{n-t}{2}\right\rfloor\) is unimodal.
\end{lemma}

The next lemma applies the generating set method to families
of sets of size at most \(k\).

\begin{lemma}[{\cite{BE00,ABEK02}}]\label{lem:prefix}
Let \(2\le t\le k\le n\) and \(2k<n+t\), and let \(r\)
be the smallest index satisfying
\[
f(r)=\max_{0\le q\le \left\lfloor\frac{n-t}{2}\right\rfloor}f(q).
\]
Then there exists a \(t\)-intersecting family
\(\calH\subseteq2^{[t+2r]}\) such that
\(\widehat{\calH}\) has maximum size among all
\(t\)-intersecting subfamilies of \(\subsetsatmost nk\).
\end{lemma}

\begin{proof}
Give an \(s\)-element set weight \(\nu_s=\mathbf 1_{\{s\le k\}}\).
Since \(f(r)\ge f(0)>0\), we have \(t+r\le k\).
In particular, \(t+2r<n\) follows from \(2k<n+t\).
By Lemma~\ref{lem:unimodal}, the standard families have
nonincreasing weights \(f(r)\ge f(r+1)\ge\cdots\ge f(\left\lfloor\frac{n-t}{2}\right\rfloor)\)
from index \(r\) onward.

For each initial segment \([u]\), the induced weight of a trace
of size \(s\) is
\[
\nu^{(u)}_s
=\sum_{j=0}^{n-u}\binom{n-u}{j}\nu_{s+j}
=B(n-u,k-s).
\]
These weights are nonnegative and nonincreasing in \(s\). In particular,
\[
\nu^{(u)}_{s-1}\nu^{(u)}_{u+t-s-1}
\ge \nu^{(u)}_s\nu^{(u)}_{u+t-s}
\qquad\bigl(t+2r+1\le u\le n,\ t\le s<(u+t)/2\bigr).
\]
This is the product condition in the generating set reduction of Bey
and Engel~\cite[Theorem~15]{BE00} (see also \cite[Proposition~5]{Bey17}).
The reduction gives a maximum-weight
family on \([n]\) whose members intersect in at least \(t\)
points within \([t+2r]\), as in~\cite[proof of Theorem~12,
Step~1]{ABEK02}.
Its traces on \([t+2r]\) form a \(t\)-intersecting family
\(\calH\). The members of size at most \(k\) of the
maximum-weight family belong to \(\widehat{\calH}\).
Since \(\widehat{\calH}\) is \(t\)-intersecting, it too
has maximum size in \(\subsetsatmost nk\).
\end{proof}

\subsection{Pushing and pulling}\label{sec:exceptional}

Ahlswede and Khachatrian~\cite{AK99} introduced the pushing--pulling
method in 1999 and used it to give new proofs of the complete
intersection theorem and Katona's intersection theorem.
Bey~\cite{Bey17} applied it to a weighted problem on two levels
of the Boolean lattice, and Filmus~\cite{Filmus17} used it in the
weighted complete intersection theorem. Wu, Li, Wu and
Feng~\cite{WLWF26} combined pushing--pulling with the generating set
method to determine the maximum total size of non-empty pairwise
cross-\(t\)-intersecting families in \(2^{[n]}\).

Let \(1\le m\le n\), and let \(\calH\subseteq 2^{[m]}\) be
shifted and \(t\)-intersecting.
For \(i,j\in[m]\), let \(A_{i,j}\) be obtained from \(A\) by
interchanging \(i\) and \(j\). 
We say that \(\calH\) is \emph{invariant under exchange operations}
on \([\ell]\) if \(A_{i,j}\in\calH\) for every
\(A\in\calH\) and \(i,j\in[\ell]\).

Suppose that \(\calH\) is not invariant on \([m]\), and let
\(\ell<m\) be the largest index for which it is invariant on
\([\ell]\). The family of \emph{exceptional members} is
\[
 \calH'=\{A\in\calH:A_{i,\ell+1}\notin\calH
                     \text{ for some }i\in[\ell]\}.
\]
It is nonempty by the choice of \(\ell\). For \(A\subseteq[m]\),
write \(A^-=A\cap[\ell]\) for its \emph{prefix} and
\(A^+=A\cap[\ell+2,m]\) for its \emph{tail}.

\begin{lemma}[{Ahlswede-Khachatrian,~\cite{AK99}}]
\label{lem:ak-exceptional}
For every \(A\in\calH'\), the following hold.
\begin{enumerate}
\renewcommand{\labelenumi}{\textup{(\roman{enumi})}}
\item \(\ell+1\notin A\), so \(A=A^-\cup A^+\).
\item \(A_{i,\ell+1}\notin\calH\) for every \(i\in A^-\).
\item \(P\cup A^+\in\calH'\) for every
      \(P\in\binom{[\ell]}{|A^-|}\).
\item If \(B\in\calH\setminus\calH'\), then
      \(|A_{i,\ell+1}\cap B|\ge t\) for every \(i\in[\ell]\).
\item If \(B\in\calH'\) and \(|A^-|+|B^-|\ne\ell+t\), then
      \(|A\cap B|\ge t+1\).
\end{enumerate}
\end{lemma}

For \(i\in[\ell]\), put
\[
 \calH'_i=\{A\in\calH':|A^-|=i\},\qquad
 \calH_i^*=\{A_{u,\ell+1}:A\in\calH'_i,\ u\in A^-\}.
\]
Take both families to be empty for \(i\notin[\ell]\).
By Lemma~\ref{lem:ak-exceptional}(i)--(iii), for each tail \(Y\)
occurring in \(\calH'_i\), the original members and their
replacements with that tail are, respectively,
\[
 \bigl\{P\cup Y:P\in\binom{[\ell]}i\bigr\},\qquad
 \bigl\{P\cup\{\ell+1\}\cup Y:P\in\binom{[\ell]}{i-1}\bigr\}.
\]
In particular, \(\calH_i^*\cap\calH=\varnothing\).

\begin{lemma}\label{lem:exceptional}
Let \(t\ge2\) and fix any positive weights \(w_0,\ldots,w_m\).
Suppose that the shifted family \(\calH\) has maximum total weight
among the \(t\)-intersecting families on \([m]\). Then every member of
\(\calH'\) has the same prefix size \(p=(\ell+t)/2\).
Moreover, for every \(A\in\calH'\) there is a \(B\in\calH'\)
such that \(B^+=[\ell+2,m]\setminus A^+\).
\end{lemma}

\begin{proof}
Since all weights are positive, maximality makes \(\calH\)
inclusion-maximal and hence upward closed.

Fix \(i\) with \(\calH'_i\ne\varnothing\), and put
\(j=\ell+t-i\). Suppose \(i\ne j\).
For \(A\in\calH'_i\) and \(u\in A^-\),
Lemma~\ref{lem:ak-exceptional}(iv),(v) gives
\[
\begin{aligned}
|A_{u,\ell+1}\cap B|
&\ge t,
&& B\in\calH\setminus\calH',\\
|A_{u,\ell+1}\cap B|
&\ge |A\cap B|-1\ge t,
&& B\in\calH'\setminus\calH'_j.
\end{aligned}
\]
Moreover, since \(2i\ne\ell+t\), for
\(A,B\in\calH'_i\), \(u\in A^-\), and \(v\in B^-\),
\[
|A_{u,\ell+1}\cap B_{v,\ell+1}|
\ge |A\cap B|-2+1
\ge t.
\]
Thus \((\calH\setminus\calH'_j)\cup\calH_i^*\)
is \(t\)-intersecting. Since
\(\calH_i^*\cap\calH=\varnothing\), maximality gives
\[
0<w(\calH_i^*)\le w(\calH'_j).
\]
In particular, \(\calH'_j\ne\varnothing\).
Interchanging \(i\) and \(j\) similarly gives
\[
w(\calH_j^*)\le w(\calH'_i).
\]

For each fixed tail, the original members and their
replacements have the same size. Counting their prefix
choices therefore gives
\[
\frac{w(\calH_i^*)}{w(\calH'_i)}
=\frac{\binom{\ell}{i-1}}{\binom{\ell}{i}}
=\frac{i}{\ell-i+1},
\qquad
\frac{w(\calH_j^*)}{w(\calH'_j)}
=\frac{j}{\ell-j+1}.
\]
Since \(i+j=\ell+t\), we obtain
\[
\begin{aligned}
1
&\ge
\frac{w(\calH_i^*)w(\calH_j^*)}
     {w(\calH'_i)w(\calH'_j)}
=\frac{ij}{(\ell-i+1)(\ell-j+1)}\\
&=
1+\frac{(\ell+1)(t-1)}
       {(\ell-i+1)(\ell-j+1)}
>1,
\end{aligned}
\]
a contradiction. Hence \(i=j\), proving
\[
|A^-|=p=\frac{\ell+t}{2}
\qquad(A\in\calH').
\]

Now fix \(A\in\calH'\) and \(u\in A^-\).
The replacement \(A_{u,\ell+1}\notin\calH\) has size
\(|A|\ge t\). By inclusion maximality and
Lemma~\ref{lem:ak-exceptional}(iv), some \(B\in\calH'\)
satisfies
\[
\begin{aligned}
t
&>|A_{u,\ell+1}\cap B|\\
&\ge (p-1)+p-\ell+|A^+\cap B^+|\\
&=t-1+|A^+\cap B^+|.
\end{aligned}
\]
Thus \(A^+\cap B^+=\varnothing\).

Put \(Y=[\ell+2,m]\setminus A^+\), so \(B^+\subseteq Y\).
Choose \(P\in\binom{[\ell]}p\) with
\[
|P\cap A^-|=2p-\ell=t.
\]
By Lemma~\ref{lem:ak-exceptional}(iii) and upward closure,
\[
P\cup B^+\in\calH',
\qquad
P\cup Y\in\calH.
\]
For any \(v\in P\cap A^-\),
\[
|(P\cup Y)_{v,\ell+1}\cap A|
=|P\cap A^-|-1
=t-1.
\]
Hence \((P\cup Y)_{v,\ell+1}\notin\calH\), so
\(P\cup Y\in\calH'\). Its tail is
\[
(P\cup Y)^+=Y=[\ell+2,m]\setminus A^+,
\]
proving the second assertion.
\end{proof}

Continue under the hypotheses of Lemma~\ref{lem:exceptional}.
For \(F\in\widehat{\calH}\) with \(F\cap[m]\in\calH'\), write
\(W=F\cap[\ell+2,n]\) for its tail in \([n]\).
By Lemmas~\ref{lem:ak-exceptional} and~\ref{lem:exceptional},
the members with exceptional trace and tail \(W\) are precisely
\[
 P\cup W,\qquad P\in\binom{[\ell]}p.
\]

Choose an intersecting family of nonempty sets
\[
 \mathcal W_1\subseteq
 \bigl\{F\cap[\ell+2,n]:F\in\widehat{\calH},
                         F\cap[m]\in\calH'\bigr\}.
\]
The remaining tails form the family
\[
 \mathcal W_0=
 \bigl\{F\cap[\ell+2,n]:F\in\widehat{\calH},
                         F\cap[m]\in\calH'\bigr\}
 \setminus\mathcal W_1.
\]
Define
\[
\begin{aligned}
\mathcal G
 &=\bigl\{P\cup\{\ell+1\}\cup W:
          P\in\binom{[\ell]}{p-1},\ W\in\mathcal W_1\bigr\},\\
\calF^*
 &=\Bigl(\widehat{\calH}\setminus
         \bigl\{P\cup W:P\in\binom{[\ell]}p,
                        W\in\mathcal W_0\bigr\}\Bigr)
   \cup\mathcal G.
\end{aligned}
\]

\begin{lemma}\label{lem:exchange}
The family \(\calF^*\subseteq\subsetsatmost nk\) is
\(t\)-intersecting, and
\[
 |\calF^*|-|\widehat{\calH}|
 =\binom{\ell}{p-1}|\mathcal W_1|
  -\binom{\ell}{p}|\mathcal W_0|.
\]
\end{lemma}

\begin{proof}
For each \(W\in\mathcal W_1\), the original and added members
with tail \(W\) have the same size \(p+|W|\le k\).
Every member of \(\mathcal G\) has trace in \(\calH_p^*\), so
\[
\mathcal G\cap\widehat{\calH}=\varnothing.
\]

Let \(A\in\mathcal G\) and
\(F\in\widehat{\calH}\cap\calF^*\).
If \(F\cap[m]\notin\calH'\), then
Lemma~\ref{lem:ak-exceptional}(iv) gives
\[
|A\cap F|\ge|A\cap F\cap[m]|\ge t.
\]
If \(F\cap[m]\in\calH'\), both tails belong to
\(\mathcal W_1\) and hence intersect. Therefore
\[
|A\cap F|
\ge(p-1)+p-\ell+1
=t.
\]
For \(A,B\in\mathcal G\), their tails also intersect,
and both sets contain \(\ell+1\). Thus
\[
|A\cap B|
\ge2(p-1)-\ell+1+1
=t.
\]
The retained members are already \(t\)-intersecting,
so \(\calF^*\) is \(t\)-intersecting.

Distinct tails give disjoint blocks of sets. Counting
the prefix choices gives
\[
|\mathcal G|
=\binom{\ell}{p-1}|\mathcal W_1|,
\qquad
|\widehat{\calH}\setminus\calF^*|
=\binom{\ell}{p}|\mathcal W_0|.
\]
Consequently,
\[
|\calF^*|-|\widehat{\calH}|
=|\mathcal G|-|\widehat{\calH}\setminus\calF^*|
=\binom{\ell}{p-1}|\mathcal W_1|
 -\binom{\ell}{p}|\mathcal W_0|.
\]
\end{proof}

\subsection{Counting and comparing extensions}

By \eqref{eq:traceweight},
\[
f(q)=\sum_{i=t+q}^{t+2q}
\binom{t+2q}{i}B(n-t-2q,k-i).
\]
We use the following adjacent comparison from~\cite{ABEK02}.

For \(1\le q\le \left\lfloor\frac{n-t}{2}\right\rfloor\), put \(a=n-t-2q\) and \(b=k-t-q\).
The sets in
\(\calS_q^{\le k}\setminus\calS_{q-1}^{\le k}\)
contain both new points \(t+2q-1,t+2q\), whereas those
in the reverse difference contain neither. Counting gives
\[
\begin{aligned}
|\calS_q^{\le k}\setminus\calS_{q-1}^{\le k}|
&=\binom{t+2q-2}{t+q-2}B(a,b),\\
|\calS_{q-1}^{\le k}\setminus\calS_q^{\le k}|
&=\binom{t+2q-2}{t+q-1}B(a,b+1).
\end{aligned}
\]
Hence
\begin{equation}\label{eq:difference}
f(q)-f(q-1)
=\binom{t+2q-2}{t+q-2}
\left(B(a,b)-\frac{q}{t+q-1}B(a,b+1)\right).
\end{equation}

For the next comparison, suppose a tail is extended by adding
at most \(d\) points from a set of \(a\ge1\) exterior points.
Fix one of these points. The numbers of extensions containing
and avoiding it are \(B(a-1,d-1)\) and \(B(a-1,d)\),
respectively. The next lemma bounds the ratio of the combined
counts for two tails whose extension limits have a fixed sum.

\begin{lemma}\label{lem:pair}
Let \(a,b,j\) be integers with \(a>2b\ge0\) and \(j\ge0\). Then
\begin{equation}\label{eq:pair}
\frac{B(a-1,b+j)+B(a-1,b-j-1)}
     {B(a-1,b+j+1)+B(a-1,b-j)}
\ge \frac{B(a,b)}{B(a,b+1)}.
\end{equation}
\end{lemma}

\begin{proof}
Both denominators are positive.
By \(B(h,s+1)=B(h,s)+\binom{h}{s+1}\),
inequality \eqref{eq:pair} is equivalent to
\[
\begin{aligned}
&\bigl(B(a-1,b+j)+B(a-1,b-j-1)\bigr)\binom{a}{b+1}\\
&\qquad\ge
B(a,b)\left(\binom{a-1}{b+j+1}+\binom{a-1}{b-j}\right).
\end{aligned}
\]

For each integer \(0\le v\le a\), expanding the binomial
coefficients gives
\[
\frac{\binom{a-1}{v+j}}{\binom{a}{v}}
=\frac{a-v}{a}
 \prod_{i=1}^{j}\max\left\{0,\frac{a}{v+i}-1\right\},
\]
where an empty product equals \(1\).
When \(v+j\ge a\), both sides vanish.
The right-hand side is convex on \([0,a]\), being a product
of nonnegative, nonincreasing convex functions. Indeed, the product
of two such functions \(g,h\) is again nonnegative and nonincreasing.
Its convexity follows from those of \(g,h\) and the inequality
\((g(x)-g(y))(h(x)-h(y))\ge 0\).

By binomial symmetry, replacing \(v\) by \(a-v\) gives
\(\binom{a-1}{v-j-1}/\binom{a}{v}\).
Thus
\[
\frac{\binom{a-1}{v+j}+\binom{a-1}{v-j-1}}{\binom{a}{v}}
\]
agrees on the integers with a convex function symmetric
about \(a/2\). Since \(a\ge2b+1\),
\[
\left|v-\frac a2\right|
\ge
\left|b+1-\frac a2\right|
\qquad(0\le v\le b).
\]
The quotient at \(v\) is therefore at least its value at
\(b+1\), giving
\[
\begin{aligned}
&\left(\binom{a-1}{v+j}+\binom{a-1}{v-j-1}\right)
 \binom{a}{b+1}\\
&\qquad\ge
\binom{a}{v}
\left(\binom{a-1}{b+j+1}+\binom{a-1}{b-j}\right).
\end{aligned}
\]
For \(v<0\), this inequality follows from nonnegativity and
\(\binom{a}{v}=0\).
Summing over \(-j\le v\le b\) gives the required inequality.
\end{proof}

The estimate can fail when \(a=2b\).
For \((a,b,j)=(2,1,1)\), the two sides of \eqref{eq:pair}
are \(2/3\) and \(3/4\), respectively.

\section{Proof of the main result (Theorem~\ref{thm:main})}\label{sec:proof}

The case \(t=1\) follows from
Proposition~\ref{prop:truncated-intersecting}.
Henceforth assume \(t\ge2\).

Let \(r\) be the least index attaining
\(\max_{0\le q\le \left\lfloor\frac{n-t}{2}\right\rfloor}f(q)\).
We show that \(\calS_r^{\le k}\) has maximum size among all
\(t\)-intersecting subfamilies of \(\subsetsatmost nk\).
Since
\[
f(r)\ge f(0)=B(n-t,k-t)>0,
\]
we have \(t+r\le k\), and hence
\[
n-t-2r>2(k-t-r)\ge0.
\]
In particular, \(n\notin[t+2r]\).

By Lemma~\ref{lem:prefix}, choose a maximum-size
\(t\)-intersecting family \(\calF\subseteq\subsetsatmost nk\)
of the form \(\calF=\widehat{\calH}\), where
\(\calH\subseteq2^{[t+2r]}\) is \(t\)-intersecting.
Among all such pairs \((\calF,\calH)\), choose one with
\(|\calH|\) largest.
If \(r=0\), then \(\calH=\{[t]\}\) and
\(\calF=\calS_0^{\le k}\), as required. Assume \(r>0\).
Shifts preserve \(t\)-intersection and the number of traces
of each size, so we may take \(\calH\) to be shifted.

Taking \(m=t+2r\) in \eqref{eq:trace-rank-weight}, we obtain
\[
w(\calH)
=\sum_{A\in\calH}w_{|A|}
=(2^{t+2r}+1)|\calF|+|\calH|.
\]
Since \(|\calH|\le2^{t+2r}\), our choice maximizes
\(w(\calH)\). All weights are positive, so \(\calH\)
is inclusion-maximal among the \(t\)-intersecting families
on \([t+2r]\).

If \(\calH\) is invariant under exchange operations on
\([t+2r]\), it contains every \(i\)-subset of \([t+2r]\)
for each size \(i\) that occurs. The minimum intersection
size of two such sets is \(\max\{0,2i-(t+2r)\}\), so
\(i\ge t+r\).
Conversely, all subsets of \([t+2r]\) of size at least
\(t+r\) together form a \(t\)-intersecting family.
Inclusion maximality therefore gives
\[
\calH=\{A\subseteq[t+2r]:|A|\ge t+r\},
\]
and consequently
\[
\calF=\widehat{\calH}=\calS_r^{\le k}.
\]

Suppose that \(\calH\) is not invariant on \([t+2r]\).
Let \(\ell<t+2r\) be the largest index for which it is
invariant on \([\ell]\), and let \(\calH'\ne\varnothing\)
be the corresponding family of exceptional traces.
By Lemma~\ref{lem:exceptional},
\[
|A^-|=p=\frac{\ell+t}{2}
\qquad(A\in\calH').
\]
Since \(p\le\ell<t+2r\) and \(p\) is an integer,
\[
t\le p\le t+r-1.
\]

Following~\cite[proof of Theorem~12, Step~2]{ABEK02}, partition the tails
of members with exceptional trace according to whether they contain \(n\):
\[
\begin{aligned}
\mathcal W_1
&=\bigl\{F\cap[\ell+2,n]:
         F\in\calF,\ F\cap[t+2r]\in\calH',\ n\in F\bigr\},\\
\mathcal W_0
&=\bigl\{F\cap[\ell+2,n]:
         F\in\calF,\ F\cap[t+2r]\in\calH',\ n\notin F\bigr\}.
\end{aligned}
\]
Every member of \(\mathcal W_1\) contains \(n\), so
\(\mathcal W_1\) is intersecting.
Lemma~\ref{lem:exchange} gives a \(t\)-intersecting family
\(\calF^*\subseteq\subsetsatmost nk\), where
\[
\begin{aligned}
\calF^*
={}&\Bigl(\calF\setminus
 \bigl\{P\cup W:
        P\in\binom{[\ell]}p,\ W\in\mathcal W_0\bigr\}\Bigr)\\
&{}\cup
 \bigl\{P\cup\{\ell+1\}\cup W:
        P\in\binom{[\ell]}{p-1},\ W\in\mathcal W_1\bigr\}.
\end{aligned}
\]

By Lemma~\ref{lem:exceptional}, the distinct tails of exceptional
traces occur in complementary pairs in \([\ell+2,t+2r]\).
This interval has positive odd size
\[
t+2r-\ell-1=2(t+r-p)-1\ge1,
\]
so the two tails in each pair are distinct.
There is at least one pair because \(\calH'\ne\varnothing\).

For each pair \(\{Y,Y'\}\), define
\[
e_i(Y,Y')
=
\bigl|\{W\in\mathcal W_i:
W\cap[\ell+2,t+2r]\in\{Y,Y'\}\}\bigr|,
\qquad i\in\{0,1\}.
\]
Put
\[
d=k-p-|Y|,
\qquad
d'=k-p-|Y'|.
\]
The exterior sets \(Z\subseteq[t+2r+1,n]\) extending
\(Y\) and \(Y'\) satisfy \(|Z|\le d\) and \(|Z|\le d'\),
respectively. Including \(n\) uses one of these places.
Counting choices from the other \(n-t-2r-1\) exterior points gives
\[
e_i(Y,Y')
=
B(n-t-2r-1,d-i)+B(n-t-2r-1,d'-i),
\qquad i\in\{0,1\}.
\]

Since \(Y,Y'\) partition \([\ell+2,t+2r]\) and
\(2p=\ell+t\),
\[
\begin{aligned}
d+d'
&=2(k-p)-(|Y|+|Y'|)\\
&=2(k-p)-(t+2r-\ell-1)\\
&=2(k-t-r)+1\ge1.
\end{aligned}
\]
Thus at least one of \(d,d'\) is nonnegative.
For the corresponding tail, \(Z=\varnothing\) gives an
extension avoiding \(n\), so
\[
e_0(Y,Y')\ge1.
\]

Interchanging \(Y,Y'\) if necessary, assume \(d\ge d'\).
Their sum then gives
\[
d=k-t-r+1+j,
\qquad
d'=k-t-r-j
\]
for some integer \(j\ge0\).
Apply Lemma~\ref{lem:pair} with
\[
a=n-t-2r,\qquad b=k-t-r.
\]
Its hypotheses hold because \(n-t-2r>2(k-t-r)\ge0\).
We obtain
\begin{equation}\label{eq:gain}
\begin{aligned}
\frac{e_1(Y,Y')}{e_0(Y,Y')}
&\ge
\frac{B(n-t-2r,k-t-r)}
     {B(n-t-2r,k-t-r+1)}\\
&>
\frac{r}{t+r-1}
\ge
\frac{p-t+1}{p}
=
\frac{\binom{\ell}{p}}{\binom{\ell}{p-1}}.
\end{aligned}
\end{equation}
The strict inequality follows from \eqref{eq:difference}
at \(q=r\), since the choice of \(r\) gives
\(f(r)>f(r-1)\). The next inequality follows from
\[
\frac{r}{t+r-1}-\frac{p-t+1}{p}
=
\frac{(t-1)(t+r-1-p)}{p(t+r-1)}
\ge0,
\]
and the last equality uses \(\ell=2p-t\).

By \eqref{eq:gain}, each complementary pair satisfies
\[
\binom{\ell}{p-1}e_1(Y,Y')
>
\binom{\ell}{p}e_0(Y,Y').
\]
Each \(W\in\mathcal W_i\) is counted for exactly one pair,
determined by its trace tail \(W\cap[\ell+2,t+2r]\).
Therefore
\[
\sum_{\{Y,Y'\}}e_i(Y,Y')=|\mathcal W_i|,
\qquad i\in\{0,1\},
\]
where the sum runs over unordered complementary pairs.
Summing the strict inequalities and applying
Lemma~\ref{lem:exchange} gives
\[
|\calF^*|-|\calF|
=
\binom{\ell}{p-1}|\mathcal W_1|
-\binom{\ell}{p}|\mathcal W_0|
>0.
\]
This contradicts the maximality of \(|\calF|\).
Hence \(\calH\) is invariant on \([t+2r]\), and
\(\calF=\calS_r^{\le k}\) attains the maximum.
\hfill\qedsymbol

The corresponding equality problem is to describe all largest
families up to permutations of \([n]\).
For \((n,k,t)=(5,3,2)\), both \(\calS_0^{\le3}\) and
\(\calS_1^{\le3}\) have size four.
The former contains a two-element set, whereas the latter
consists entirely of three-element sets, so they are not
isomorphic.

\paragraph{Declaration of generative AI use.}
The work was carried out with the aid of generative AI tools. The authors devised the rough proof strategies, while the generative AI tools filled in most of the computational details; the resulting content was then checked and rewritten by the authors.

\end{document}